\documentclass{article}
\usepackage{mathtext}
\usepackage[T2A]{fontenc}
\usepackage[cp866]{inputenc}
\usepackage[russian,english]{babel}
\usepackage{amssymb}
\usepackage{amsmath,amsthm}
\usepackage{titlesec}
\usepackage{titletoc}

\begin{document}

\textwidth 6.6in \textheight 8.6in \footskip 0.3in
\parskip 0.06in

\renewcommand{\baselinestretch}{1.0}

\oddsidemargin 0in \evensidemargin 0in

\baselineskip 16pt

\title{\bf Local definitions of formations of finite groups}

\author{L.A.  Shemetkov\\ 
{\small Department of Mathematics, F.~Scorina Gomel State University,} \\
{\small Gomel 246019, Belarus;} \\
{\small E-mail: shemetkov@gsu.by}}

\date{}
\maketitle

\begin{abstract}

A problem of constructing of local definitions for formations of finite groups is
discussed in the article. The author analyzes relations between local 
definitions of various types. A new proof of existence  of an
$\omega$-composition satellite of an $\omega$-solubly saturated formation is 
obtained. It is proved that if a non-empty formation of finite groups is
$\mathfrak X$-local by F\"{o}rster, then  it has an $\mathfrak X$-composition 
satellite.

\end{abstract}

{\bf 2000 Mathematics Subject Classification:} 20D10

{\bf Keywords:} Frattini subgroup, saturated formation, satellite

\bigskip
 \textbf{1. Introduction }

\bigskip 
We consider only finite groups. So, all group classes considered are subclasses
of the class  $\mathfrak E$ of all finite groups.
 Recall that a formation is a group class closed under taking homomorphic 
images and subdirect  products (see \cite{Gas}). A formation  $\mathfrak F$
 is said to be  $p$-saturated ($p$ a prime) if the condition:
\[
 G/N\in \mathfrak F\; \text  {for a}\; G\text{-invariant}\;  p\text 
{-subgroup}\; N\; \text {of}\; \Phi (G)
\]
 always implies
$G\in \mathfrak F$. A formation $\mathfrak F$
 is said to be ${\mathfrak N}_p$-saturated if the condition
\[
G/\Phi (N)\in \mathfrak F \; \text {for a}\; \text {normal}\,  p\text 
{-subgroup}\, N\, \text  {of}\; G
\]
 always implies $G\in \mathfrak F$.

If a formation is  $p$-saturated for any prime  
 $p$, then it is called saturated. 
Clearly, every  $p$-saturated formation is  ${\mathfrak N}_p$-saturated.
 The converse is not true: there is  an extensive class of
${\mathfrak N}_p$-saturated formations which are not  $p$-saturated.
However, as it is established in \cite{Shem97},
 between local definitions of these two types of formations there is a 
close connection.

The concept  of local definitions
 of saturated formations was 
considered for the first time by W. Gasch\"utz \cite{Gas}.
 Following \cite{Shem2001}, we
 formulate it in the 
general form.

A local definition is a map
$f$:$\mathfrak E$ $\to $\{formations\}
together with a $f$-rule which decide whether a chief factor is 
$f$-central or $f$-eccentric in a group. In addition, we follow the 
agreement that the local definition $f$ does not distinguish between non-identity groups with 
the same (up to isomorphism)  set of composition factors. 
Therefore, for any fixed prime $p$, $f$ is not distinguish between any two non-identity 
$p$-groups;  we will denote through $f(p)$ a value of $f$ on non-identity 
$p$-groups.

If a class $\mathfrak F$ coincides with the class of all groups all of whose  
chief factors are $f$-central, we say that $f$ is a local definition of 
$\mathfrak F$.
It generalises the  concept of nilpotency. Thus, the problem of  finding 
 local definitions for group classes is equivalent to a problem of 
 finding classes  of generalized nilpotent groups.

In this paper we analyze relations between local definitions of different 
types and give a new proof of a theorem on a local definition of a 
formation
which is ${\mathfrak N}_p$-saturated for any  $p$ in a set $\omega$ of primes.

\bigskip
\textbf{2. Preliminaries}

\bigskip
We use standard notations and definitions \cite{DH}. We say that a map $f$ does not distinguish between  $\mathfrak H$-groups if
 $f(A)=f(B)$ for any two
groups $A$ and $B$ in $\mathfrak H$. 
 Following Gasch\"{u}tz, the $\mathfrak F$-residual
$G^{\mathfrak F}$  of a group $G$ is the least normal subgroup with
 quotient  in $\mathfrak F$.
 The Gasch\"{u}tz product $\mathfrak F \circ \mathfrak H$ of formations
 $\mathfrak F$
and $\mathfrak H$ is defined as the class of all groups
 $G$ such that $G^{\mathfrak H} \in \mathfrak F$.
If $\mathfrak F$ is closed under taking of normal subgroups, then
$\mathfrak F \circ \mathfrak H$ coincides with the class
 $\mathfrak F  \mathfrak H$ of all extensions of  $\mathfrak F$-groups by
$\mathfrak H$-groups.

$\mathbb{P}$ is the set of all primes;
${\mathrm {Char}}(\mathfrak X)$ is the set of orders of all simple abelian
 groups 
in $\mathfrak X$. 
A group $G$ is called a $pd$-group if its order is divisible by a
prime $p$; $C_p$ is a group of order $p$; if $\omega \subseteq 
\mathbb{P}$, then $\omega '= \mathbb{P}\setminus \omega$;
 an $\omega d$-group (a chief $\omega d$-factor) is a group 
(a chief factor) being $pd$-group for some $p \in \omega$;  $G_{\omega d}$
 is the largest normal subgroup all of whose
 $G$-chief  factors are $\omega d$-groups ($G_{\omega d}=1$ if all minimal
 normal subgroups in $G$ are $\omega '$-groups). If  $\mathfrak H$ is a class
 of groups, then ${\mathfrak H}_\omega$ is the class of all $\omega$-groups in 
$\mathfrak H$.  A chief factor $H/K$ of $G$ is called a chief $\mathfrak 
H$-factor if $H/K\in \mathfrak H$. The socle $\mathrm{Soc}(G)$ of a group
$G\ne 1$ is the product of all minimal normal subgroups of $G$.

  $[A]B$ is a semidirect product with a normal subgroup $A$;
  $O_{\omega}(G)$ is the largest normal $\omega$-subgroup in
  $G$; $\pi (G)$ is the set of all primes dividing the order of a group $G$;
 $\pi (\mathfrak{F)}=\cup _{G\in \mathfrak{F}}\pi (G)$;
 $\mathfrak{N}$ is the class of all nilpotent groups; $\mathfrak{A}$ is the class
of all abelian groups;  $\mathrm{Com}(G)$ is the class of all groups
 that are isomorphic to composition factors of a group $G$;
 $\mathrm{Com}(\mathfrak F)=\cup _{G\in \mathfrak{F}}\mathrm{Com}(G)$;
 $\mathrm{Com}^+(\mathfrak F)$ is the class of all abelian groups in 
$\mathrm{Com}(\mathfrak F)$; $\mathrm{Com}^-(\mathfrak F)$ is the class of all non-abelian 
groups in $\mathrm{Com}(\mathfrak F)$;
 $(G)$ is the class of all groups isomorphic to $G$; $\mathfrak{J}$ is the class
 of all simple (abelian and non-abelian) groups; if $\mathfrak{L}$ is a subclass
 in $\mathfrak{J}$, then    $\mathfrak L '= \mathfrak J \setminus \mathfrak L$;
 $\mathfrak{L}^{+}$ is the class of all abelian groups in $\mathfrak{L}$,
 $\mathfrak{L}^{-}=\mathfrak{L}\setminus \mathfrak{L}^{+}$.
$\mathrm{E}\mathfrak H$ is  the class of all groups $G$ such that
$\mathrm{Com}(G)\subseteq \mathfrak H$; $G_{\mathrm{E}\mathfrak H}$
is the $\mathrm{E}\mathfrak H$-radical of $G$, the largest normal
$\mathrm{E}\mathfrak H$-subgroup in $G$. If $S\in \mathfrak J$, then $C^S(G)$
is the intersection  of centralizers of all chief $\mathrm{E}(S)$-factors 
of $G$ ($C^S(G)=G$ if $S\notin \mathrm{Com}(G)$); if $S=C_p$, we write
$C^p(G)$ in place of $C^S(G)$.

\medskip
\textbf{Lemma 2.1} (see \cite{Shem2001}, Lemmas 2--3).
\textsl{$(a)$ If $S$ is a non-abelian simple group, then $C^S(G)$ is the 
$\mathrm{E}(S)'$-radical of $G$, the largest normal subgroup not having 
composition factors isomorphic to $S$.}

\textsl{$(b)$ Let $p$ be a prime, and $\mathfrak H$  be  the class of all
 groups all of whose chief
$p$-factors are central. Then $C^p(G)$ is the $G_{\mathfrak H}$-radical of $G$,
for every group $G$.}  
   
\medskip
The following three lemmas are reformulations of Lemmas IV.4.14--IV.4.16 
in \cite{DH}
whose proofs use only $p$-solubly saturation.

\medskip
\textbf{Lemma 2.2.} \textsl{Let $\mathfrak F$ be an ${\mathfrak N}_p$-saturated 
formation, $p$ a prime. If $C_p\in \mathrm{Com}(\mathfrak F)$, then
${\mathfrak N}_p\subseteq \mathfrak F$.}

\medskip
\textbf{Lemma 2.3.} \textsl{Let $\mathfrak F$ be an ${\mathfrak N}_p$-saturated 
formation containing ${\mathfrak N}_p$, $p$ a prime.  Let $N$ be an 
elementary abelian normal $p$-subgroup in $G$ such that 
$[N](G/N)\in \mathfrak F$. Then $G\in \mathfrak F$.}

\medskip
\textbf{Lemma 2.4.} \textsl{Let $p$ be a prime, and let $\mathfrak F$ be 
an ${\mathfrak N}_p$-saturated 
formation containing ${\mathfrak N}_p$. Let $N$ be an elementary abelian normal
 $p$-subgroup in $G$ such that  $G/N\in \mathfrak F$
 and $[N](G/C_G(N))\in \mathfrak F$. Then $G\in \mathfrak F$.}

\medskip
\textbf{Proof.} Set  $M=[N](G/N)$, $C=C_G(N)$. Evidently, $C/N=C_{G/N}(N)$.
In the group $M$ we have $C_M(N)=N\times C/N$ and $C/N$ is normal in $M$. Hence
$M/(C/N)\simeq [N](G/C)\in \mathfrak F$. Since $M/N\in \mathfrak F$, it 
follows that $M/N\cap(C/N)\simeq M\in \mathfrak F$. Now we apply Lemma 
2.3. \qed

\medskip
\textbf{Lemma 2.5.} \textsl{Let $\mathfrak F$ be an ${\mathfrak N}_p$-saturated 
formation containing ${\mathfrak N}_p$, $p$ a prime. Let $H\in \mathfrak F$ and let
$C^p(H)\le L \unlhd H$. If $N$ is an irreducible ${\mathbb F}_p(H/L)$-module, then
$[N](H/L)\in \mathfrak F$. }

\medskip
\medskip
\textbf{Lemma 2.6} (see \cite{DH}, Proposition IV.1.5). \textsl{Let $\mathfrak F$ be a formation and $G\in 
\mathfrak F$. Let $S,R,K$ be normal subgroups in $G$ such that $S\subseteq R$ and
$K\subseteq C_G(R/S)$.  Then $[R/S](G/K)\in \mathfrak F$.}

\medskip
\textbf{Lemma 2.7} (see \cite{Sch72} or \cite{For}, Theorem 7.11).
\textsl{If $H/\Phi (G)=\mathrm{Soc}(G/\Phi (G)$, then $C_G(H)\subseteq H$.}

\medskip
\textbf{Lemma 2.8} (see \cite{DH}, Lemma IV.4.11).
\textsl{Let $p$ be a prime, $L=\Phi (O_p(G))$. Then $C^p(G/L)=C^p(G)/L$.}

 \bigskip
\textbf{3. Local and $\omega$-local satellites}

\bigskip
The following type of a local definition was proposed by W.~Gasch\"utz \cite{Gas}. 

\medskip
\textbf{ Definition 3.1.}
Let $f$ be a local definition such that
 \[
 f(A)=\bigcap_{p\in \pi (A)}f(p)
\]
 for any group 
$A\ne 1$. Let an $f$-rule be defined as follows: a chief 
factor $H/K$ of a group $G$ is $f$-central if $G/C_G(H/K)\in 
f(H/K)$. Then $f$ is called {\it a local satellite}.

\medskip
\textbf{ Definition 3.2} (see \cite{DH}, p. 387). Let $A$ be a group of 
operators for a group $G$, and $f$ a local satellite.

(i) We say that $A$ acts $f$-\emph{centrally} on an $A$-composition 
factor $H/K$ of $G$ if $A/C_A(H/K)\in f(p)$ for every prime $p\in \pi (H/K)$.

(ii) We say that $A$ acts $f$-\emph{hypercentrally} on $G$ if $A$ acts 
$f$-centrally on every $A$-composition factor  of $G$.

\medskip
The convenient notation $LF(f)$ for a group class with a local satellite 
$f$ was introduced by Doerk and Hawkes \cite{DH}.  Clearly, $LF(f)$ is a 
non-empty formation (we have always $ 1\in  LF(f)$).

The following proposition is evident.       

\medskip                                            
\textbf{Proposition 3.1.} \textsl{  
Let $f$ be a local satellite and ${\pi =\{p\in \mathbb{P}\ \mid \ f(p)\ne 
\varnothing\}}$. Then $LF(f)$ consists precisely of $\pi$-groups $G$ 
satisfying the following condition: $G/O_{p', p}(G)\in f(p)$ for any
 $p\in \pi (G)$. Thus, if $\pi =\varnothing$, we have $LF(f)=(1)$.  If $\pi \ne 
\varnothing$, we have that
\[
 LF(f)={\mathfrak E}_{\pi }\bigcap (\bigcap_{p\in 
\pi } ({\mathfrak E}_{p'} {\mathfrak E}_{p} f(p))).
\]
}

\medskip                        
We remember the reader that a formation $\mathfrak F$ is  saturated if 
$G/\Phi(G)\in \mathfrak F$ always implies $G\in\mathfrak F$ (by definition, 
the empty set is a saturated formation). W.~Gash\"{u}tz has shown that every formation
with a local satellite is saturated. This fact follows also from the following
 theorem of P.~Schmid.

\medskip 
\textbf{Theorem 3.1 } (see \cite{DH}, Theorem IV.6.7). \textsl{Let $f$ be a
 local satellite, and let $A$ be
a group of operators for a group $G$. If  $A$ acts $f$-hypercentrally on
 $G/\Phi(G)$, then $A$ acts likewise on $G$.}

 \medskip 
 The following remarkable result is known as the  
Gasch\"utz--Lubeseder--Schmid theorem, see \cite{DH}, Theorem 
IV.4.6.

\medskip
\textbf{Theorem 3.2. }  \textsl{  
A non-empty formation has a local satellite if and only if it is 
saturated. }
                 
\medskip
It is straightforward to verify that if $\mathfrak F$ is a non-empty formation, then 
$\mathfrak N \mathfrak F$ is a formation with a local satellite $f$ such 
that $f(p)=\mathfrak F$ for every prime $p$.
Evidently, the formation ${\mathfrak A}_p \times {\mathfrak N}_{p'}$ of all nilpotent
 groups with an abelian Sylow 
$p$-subgroup is not saturated, but for every prime $q\ne p$,
 $G/(\Phi(G)\cap O_q(G))\in {\mathfrak A}_p \times {\mathfrak N}_{p'}$ always implies $G\in {\mathfrak A}_p \times 
{\mathfrak N}_{p'}$. One more fact of the same sort is the following. 
Consider a saturated formation of the form $\mathfrak M \circ \mathfrak H$. Here 
$\mathfrak H$ can be non-saturated, but for every prime
 $p\in \mathbb P \setminus \pi (\mathfrak M)$,
$G/(\Phi(G)\cap O_p(G))\in 
\mathfrak H$ always implies $G\in \mathfrak H$. The facts of such kind 
lead to the concept of a $\omega$-saturated formation \cite{ShSk95}.

\medskip
\textbf{Definition 3.3. }
 Let $\omega$ be a set of primes. A formation $\mathfrak F$ 
is called  $\omega$-{\itshape saturated} if for every prime $p\in
 \omega$, 
$G/(\Phi(G)\cap O_p(G))\in \mathfrak F$ always implies $G\in \mathfrak F$.

\medskip
The problem of 
finding of local definitions of $\omega$-saturated formations was 
considered in \cite{ShSk99} and \cite{Shem2001}.  While solving this problem the following concept of 
small centralizer was useful (see \cite{ShBall}).

\medskip
\textbf{Definition 3.4. }
 Let $H/K$ be a chief factor of a group $G$. The
\textsl { small centralizer} $c_G(H/K)$ 
of $H/K$ in $G$ is the subgroup generated by all normal subgroups $N$ of 
$G$ such that $\mathrm{Com} (NK/K)\cap \mathrm{Com} (H/K)=\varnothing$.

\medskip
With the help of Definition 3.4 we can introduce  the concept `$\omega $-saturated satellite'
 as follows.

\medskip
\textbf{Definition 3.5. }
 Let $\omega $ be a set of primes, and $f$ a local definition which does 
not distinguish between all non-identity $\omega '$-groups; if $\omega '\ne \varnothing $, 
we denote through $f(\omega ')$ a value of $f$ on non-identity $\omega 
'$-groups. In addition, we assume that
\[
 {f(A)=\bigcap _{p\in \pi (A)\cap \omega }f(p)}
\]
 for any $\omega d$-group $A$. Let an 
$f$-rule be defined by the following way: a chief factor $H/K$ of $G$ is 
$f$-central in $G$ if either $H/K$ is an $\omega d$-group and 
$G/C_G(H/K)\in f(H/K)$ or else $H/K$ is an $\omega '$-group and 
$G/c_G(H/K)\in f(\omega ')$. Then $f$ is called an $\omega$-{\itshape local 
satellite}. We denote by $LF_{\omega }(f)$ the class of all groups all of whose 
 chief factors are $f$-central. By definition, $1\in LF_{\omega }(f)$.

\medskip
Clearly, if $\omega =\mathbb P$, then an $\omega$-local 
satellite $f$ is a local satellite,  and 
$LF_{\omega }(f)=LF(f)$. If  $\omega \ne \mathbb P$  and $f(\omega 
')=\varnothing$, then $LF_{\omega }(f)=LF(h)$   where $h(p)=f(p)$ if 
$p\in \omega$, and $h(p)=\varnothing$ if $p\in \omega '$.

\medskip
\textbf{Lemma 3.1} (see \cite{Shem2001}, Lemma 1).\textsl{  
 Let $\mathfrak L$  be a subclass in $\mathfrak J$, and
 ${\{S_i\ |\ i\in I\}}$ be the set of all 
 $\mathrm{E} \mathfrak L$-factors of a group $G$. Then
 $\bigcap _{i\in I}c_G(S_i)$  is the
  $\mathrm{E} (\mathfrak L')$-radical $G_{\mathrm{E} (\mathfrak L')}$ of $G$.}

\medskip
\textbf{Remark 3.1.} 
In Lemma 3.1 the set $\{c_G(S_i) \mid i\in I\}$ can be empty.
We always follow the agreement that the intersection of an empty set of 
subgroups of $G$ coincides with $G$.

\medskip
The following proposition is similar to Proposition 3.1.

\medskip
\textbf{Proposition 3.2.} \textsl{  
Let $f$ be an $\omega$-local satellite, and $\omega$ a proper subset in~$\mathbb P$. 
Let $\pi =\{p\in \omega \ |\ f(p)\ne \varnothing \}$. Then: }

\textsl{  
$(1)$  if $\pi =\varnothing$ and $f(\omega ')=\varnothing$, then
 $LF_{\omega }(f)=(1)$; }

\textsl{  
$(2)$ if $\pi =\varnothing$ and $f(\omega ')\ne \varnothing$, then $LF_{\omega 
}(f)={\mathfrak E}_{\omega '} \cap f(\omega ')$;}

\textsl{  
$(3)$ if $f(\omega ')\ne \varnothing$, then $LF_{\omega 
}(f)$ consists precisely of groups $G$ such that $G/G_{\omega d}\in f(\omega ')$
and $G/O_{p',p}(G)\in f(p)$ for any $p\in \pi (G)\cap \omega $.}

\medskip
\textbf{Proof.}
Statements (1) and (2) are evident.

 Prove (3). Assume that $f(\omega ')\ne \varnothing$, and
let $G\in LF_{\omega }(f)$. Let $\mathfrak T$ be the set of all chief $\omega 
'$-factors in $G$. If  a chief factor $H/K$ of $G$ is an $\omega '$-group, then 
$G/c_G(H/K)\in f(\omega ')$. Therefore,
 $G/\bigcap_{H/K \in \mathfrak T} c_G(H/K)\in f(\omega ')$. By Lemma 3.1,
 $\bigcap_{H/K \in \mathfrak T}c_G(H/K)=G_{\omega d}$. So, 
$G/G_{\omega d}\in f(\omega ')$.   If $p \in \omega$ and $H/K$ is an chief
 $pd$-factor, then
  $G/C_G(H/K)\in f(p)$, and we have   $G/O_{p',p}(G)\in f(p)$.

Conversely, let $G$ be a group such that $G/G_{\omega d}\in f(\omega ')$ 
and ${G/O_{p',p}(G)\in f(p)}$ for any $p\in \pi (G)\cap \omega$. Clearly, we 
have that all $G$-chief $\omega d$-factors are $f$-central. Let $H/K$ be a 
$G$-chief $\omega '$-factor of $G$. Then $G_{\omega d}K/K\subseteq 
c_G(H/K)$, and $G/G_{\omega d}\in f(\omega ')$ implies $G/c_G(H/K)\in 
f(\omega ')$.
\qed

\medskip
The following result extends Theorem 3.2 to $\omega$-saturated formations.

\medskip
\textbf{Theorem 3.3} (see \cite{ShSk99}, Theorem 1).\textsl{  
 Let $\omega$ be a set of primes. 
A non-empty formation has a $\omega$-local satellite if and only if it is 
$\omega$-saturated.}

\bigskip
\textbf{4. Composition and $\mathfrak L$-composition satellites }

 \bigskip
Gasch\"utz's main idea \cite{Gas} was to study groups modulo $p$-groups, and he 
implemented it through local satellites of soluble formations.
 While considering non-soluble 
formations, we have to follow the following principle: study groups modulo
$p$-groups and simple groups. That approach was proposed in the lecture
 \cite{Shem73} at the
conference in 1973; in that lecture 
  composition satellites were considered  under the name
`primarily homogeneous screens'. 

\medskip
\textbf{Definition 4.1. }
Let $f$ be a local definition, and let an $f$-rule be defined as follows:
a chief factor $H/K$ of a group $G$ is $f$-central if $G/C_G(H/K)\in 
f(H/K).$  Then $f$ is called \textsl{ a composition satellite}. We denote by
$CF(f)$ the class of all groups all of whose chief factors are $f$-central.

\medskip
\textbf{Definition 4.2.} Let $A$ be a group of 
operators for a group $G$, and $f$ a composition satellite.

(i) We say that $A$ acts $f$-\emph{centrally} on an $A$-composition 
factor $H/K$ of $G$ if $A/C_A(H/K)\in f(H/K)$.

(ii) We say that $A$ acts $f$-\emph{hypercentrally} on $G$ if it acts 
$f$-centrally on every $A$-composition factor  of $G$.

\medskip
As an example, we consider the class $\mathfrak N  ^*$ of all quasinilpotent
groups (for the definition of a quasinilpotent group,
 see \cite{Hup}, Definition X.13.2). It is easy to check that
$\mathfrak N  ^*= CF(f) $  where $f$ is a composition satellite such that
$f(p)=(1)$ for every prime $p$, and $f(S)=\mathrm{form}(S)$ for every non-abelian simple group
$S$. Here  $\mathrm{form}(S)$ is a least formation containing $S$; it consists
 of all groups represented as a direct product
 $A_1\times \dots \times A_n$ with $A_i \simeq S$ for any $i$. 
The formation $\mathfrak N  ^*$ is non-saturated, but it is solubly saturated.   

\medskip
As pointed out in \cite{DH}, formations with composition satellites
 were also con\-si\-dered---in different terminology---by R.~Baer in his
 unpublished manuscript. By R.~Baer, a formation $\mathfrak F$ is 
called  {\it solubly saturated} if the condition 
$G/\Phi (G_{\mathfrak S})\in \mathfrak F$ always implies $G\in \mathfrak 
F$ (here $G_{\mathfrak S}$ is the soluble radical of $G$).
 The question of the coincidense of the
family of non-empty solubly saturated formations and the family of 
formations with
composition satellites was solved by the following result due to R.~Baer.

\medskip
\textbf{Theorem 4.1 }(see \cite{DH}, Theorem IV.4.17). \textsl{  
A non-empty formation has a composition satellite if and only if it is 
solubly saturated.}
           
 \medskip
 A composition satellite $h$ is called  \emph{integrated} if $h(S)\subseteq CF(h)$ for
any simple group $S$. If $\mathfrak F = CF(f)$, then  $\mathfrak F = CF(h)$ 
where $h(S)=f(S) \cap \mathfrak F$ for any simple group $S$. Thus, if a formation
has a composition satellite, then it has an integrated composition satellite.

\medskip
\textbf{Remark 4.1.} Let  $\{CF(f_i)\mid i\in I\}$ be a family of 
formations having composition satellites. Let $f= \cap _{i\in I}f_i$ be a composition 
satellite such that $f(S)= \cap _{i\in I}f_i(S)$ for every
 $S\in \mathfrak J$. Clearly, $CF(f)=\cap _{i\in I}CF(f_i)$.

\medskip
\textbf{Remark 4.2.} Let $\mathfrak X$ be a set of groups. Let 
$\{{\mathfrak F}_i \mid i\in I\}$  be the class of all formations 
${\mathfrak F}_i$
satisfying the following two conditions:
 1) $\mathfrak X  \subseteq  {\mathfrak F}_i$;
 2) ${\mathfrak F}_i$ has a composition satellite. Set 
 $\mathrm{cform}(\mathfrak X)= \cap _{i\in I}{\mathfrak F}_i$.
By Remark 4.1, $\mathrm{cform}(\mathfrak X)$ has a
composition satellite. In the subsequent we will use that notation
$\mathrm{cform}(\mathfrak X)$.  

\medskip
\textbf{Remark 4.3.} Assume that a non-empty formation $\mathfrak F$
 has an composition satellite. Let $\{f_i\mid i\in I\}$ be the class of all
composition satellites of $\mathfrak F$. Having in mind Remarks 4.1 and  4.2 we see
that $f= \cap _{i\in I}f_i$ is a composition satellite of $\mathfrak F$;
$f$ is called the \emph{minimal composition satellite} of $\mathfrak F$.

\medskip
\textbf{Lemma 4.1.}\textsl{ Let $\mathfrak X$ be a set of groups, and 
$S$ a simple group. Then 
$\mathfrak H = \mathrm{Q}(G/C^S(G) \mid G\in \mathrm{form}(\mathfrak X))$ is a
 formation, and $\mathrm{Com}(\mathfrak H)\subseteq \mathrm{Com}(\mathfrak X).$ }

\medskip
\textbf{Proof.} By Proposition IV.1.10 in \cite{DH}, $\mathfrak H$ is a 
formation.  By Lemma II.1.18 in \cite{DH}, $\mathrm{form}(\mathfrak X)=
\mathrm{QR_0}\mathfrak X$.  Therefore, inclusion 
$\mathrm{Com}(\mathfrak H)\subseteq \mathrm{Com}(\mathfrak X)$ is valid. \qed

\medskip
\textbf{Lemma 4.2.}\textsl{ Let $\mathfrak X$ be a non-empty set of groups, and 
$f$ be a composition satellite such that
$f(S) = \mathrm{Q}(G/C^S(G) \mid G\in \mathrm{form}(\mathfrak X))$ if 
$S\in \mathrm{Com}(\mathfrak X)$, and $f(S)=\varnothing$
 if $S\in \mathfrak J \setminus \mathrm{Com}(\mathfrak X)$. Then
$f$ is the minimal composition satellite of $\mathrm{cform}(\mathfrak X)$. }

\medskip
\textbf{Proof.} Let $f_1$ be the minimal composition satellite of
$\mathfrak F = \mathrm{cform}(\mathfrak X)$ (see Remark 4.3). We will
prove that $f_1=f$. 

Since $\mathfrak X \subseteq \mathfrak F$, $G/C^S(G)\in f_1(S)$ for any
group $G\in \mathfrak X$ and any $S\in \mathrm{Com}(G)$ and therefore
$f(S)\subseteq f_1(S)$. So $CF(f)\subseteq \mathfrak F \subseteq CF(f_1)$.
On the other hand, $\mathfrak X \subseteq CF(f)$. Thus $\mathfrak F=CF(f)$ 
and $f=f_1$.
 \qed

\medskip
The following theorem proved independently in \cite{Grad}  and
 \cite{Sch} was the first important
result on composition formations.

\medskip
\textbf{Theorem 4.2.} \textsl{  
Let $f$ be an integrated composition satellite. Let $A$ be a group of 
automorphisms of a group $G$. If $A$ acts $f$-hypercentrally 
on $G$, then $A\in CF(f)$. }

\medskip
Applying Theorem 4.2 to the formation $\mathfrak U$ of all supersoluble groups,
we have the following result.

\medskip
\textbf{Theorem 4.3 } (see \cite{Grad}, Theorem 2.4).\textsl{  
Let $A$ be a group of automorphisms of a group $G$. Assume that there 
exists a chain of  $A$-invariant  subgroups
\[
 G=G_0 > G_1 > \dots > G_n=1
\]
with prime indices  $|G_{i-1}:G_i|$. Then $A$ is supersoluble.  }

\medskip
In 1968 S.A. Syskin tried to prove Theorem 4.3 in the soluble universe,
 but his proof \cite{Sys}  is false.

In \cite{Shem97} there has been begun studying of local definitions of
 $\omega$-solubly saturated formations.

\medskip
\textbf{Definition 4.3. } Let $\omega$ be a set of primes. 
 A formation $\mathfrak F$ is called:

(1) $\omega$-{\itshape solubly saturated} if the condition
\[
G/N\in \mathfrak F \; \text {for}\, G\text {-invariant}\,  \omega \text 
{-subgroup}\, N\, \text  {in}\; \Phi(G_{\omega{\text -}\mathfrak S})
\]
 always implies $G\in \mathfrak F$ (here $G_{\omega{\text -}\mathfrak S}$ is 
the $\omega$-soluble radical of $G$);

(2) ${\mathfrak N}_{\omega}$-{\itshape saturated} if for every prime $p\in \omega$,
the condition  $G/\Phi (O_p(G))\in \mathfrak F$
always implies $G\in \mathfrak F$.

\medskip
Later we  will establish that  the $p$-solubly saturation is equivalent to
the ${\mathfrak N}_p$-saturation, and therefore a formation $\mathfrak F$
is $\omega$-solubly saturated if and only if it is $p$-solubly saturated for
every $p\in \omega$. 

\medskip
\textbf{Definition 4.4.} Let $\mathfrak L$ be a class of simple groups. Let $f$ 
be a local definition which does  not distinguish between all non-identity 
$\mathrm{E}(\mathfrak L ')$-groups; if $\mathfrak L ' \ne \varnothing$, we 
denote by $f(\mathfrak L ')$ an value of $f$ on non-identity 
$\mathrm{E}(\mathfrak L ')$-groups. Let $f$-rule be defined as follows:
a chief factor  $H/K$ of a group $G$ is $f$-central in $G$ if either $H/K$ 
is an $\mathrm{E}\mathfrak L$-group and $G/C_G(H/K)\in f(H/K)$ or $H/K$ 
is a  $\mathrm{E}(\mathfrak L ')$-group and $G/c_G(H/K)\in f(H/K)=f(\mathfrak L ')$.  
Then $f$ is called an $\mathfrak L${\itshape -composition satellite}. We denote by
$CF_{\mathfrak L}(f)$ the class of all groups all of whose chief factors 
are $f$-central. By definition, $1\in CF_{\mathfrak L }(f)$.

\medskip
Clearly, if $\mathfrak L = \mathfrak J$, then an $\mathfrak L$-composition
 satellite  $f$ is a composition satellite,  and 
$CF_{\mathfrak L}(f)=CF(f)$. If $\mathfrak L \ne \mathfrak J$   and
 $f(\mathfrak L ')=\varnothing$, then $CF_{\mathfrak L}(f)=CF(h)$ 
  where $h(S)=f(S)$ if 
$S\in \mathfrak L$, and $h(S)=\varnothing$ if $S\in \mathfrak L '$.

\medskip
\textbf{Proposition 4.1.}\textsl{  
Let $\mathfrak L$ be a class of simple groups, and $f$ an
 $\mathfrak L$-composition satellite. Let 
 $\mathfrak K=\{S\in \mathfrak L \mid f(S) \ne \varnothing\}.$  Then:}

\textsl{  
$(1)$ if $\mathfrak K = \varnothing$ and  $f(\mathfrak L ')=\varnothing$, 
then $CF_{\mathfrak L}(f)=(1)$; }

\textsl{  
$(2)$ if $\mathfrak K = \varnothing$
 and  $f(\mathfrak L ')\ne \varnothing$, then
 $CF_{\mathfrak L}(f)= \mathrm{E}(\mathfrak L ')\cap f(\mathfrak L ')$; }

\textsl{  
$(3)$ if  $f(\mathfrak L ')\ne \varnothing$,
then $CF_{\mathfrak L}(f)$ consists precisely of groups $G$ such that
$G/G_{\mathrm{E}\mathfrak L}\in f(\mathfrak L ')$ and $G/C^S(G)\in f(S)$ 
for every $S\in \mathrm{Com}(G) \cap \mathfrak L$. }

\medskip
\textbf{Proof.}
Statements (1) and (2) are evident.

Prove (3). Assume that $f(\mathfrak L ')\ne \varnothing$, and
let $G\in CF_{\mathfrak L }(f)$. Let $\mathfrak T$ be the set of all chief
 $\mathrm{E}(\mathfrak L')$-factors in $G$. If  a chief factor $H/K$ of $G$ is
an $\mathrm{E}(\mathfrak L')$-group, then 
$G/c_G(H/K)\in f(\mathfrak L ')$. Therefore,
 $G/\bigcap_{H/K \in \mathfrak T} c_G(H/K)\in f(\mathfrak L ')$. By Lemma 3.1,
 $\bigcap_{H/K \in \mathfrak T}c_G(H/K)=G_{\mathrm{E}\mathfrak L}$. So, 
$G/G_{\mathrm{E}(\mathfrak L)}\in f(\mathfrak L ')$. 
  If $S \in \mathfrak L$ and $H/K$ is an chief
 $\mathrm{E}(\mathfrak L)$-factor, then
  $G/C_G(H/K)\in f(S)$, and we have  $G/C^S(G)\in f(S)$.

Conversely, let $G$ be a group such that 
$G/G_{\mathrm{E}\mathfrak L}\in f(\mathfrak L ')$ and 
$G/C^S(G)\in f(S)$ 
for every $S\in \mathrm{Com}(G) \cap \mathfrak L$. Clearly,  
 all chief $\mathrm{E}\mathfrak L$-factors of $G$ are $f$-central.
 Let $H/K$ be a 
chief $\mathrm{E}(\mathfrak L')$-factor of $G$. Then
 $G_{\mathrm{E}\mathfrak L}\subseteq 
c_G(H/K)$, and therefore $G/G_{\mathrm{E}\mathfrak L}\in f(\mathfrak L ')$ implies
 $G/c_G(H/K)\in f(\mathfrak L ')$. 
\qed

\medskip
An $\mathfrak L$-composition satellite $f$ is called \emph{integrated} if
$f(S)\in CF_{\mathfrak L}(f)$ for every $S\in \mathfrak J$.
 If $\mathfrak F=CF_{\mathfrak L}(f)$, then $\mathfrak F=CF_{\mathfrak L}(h)$ 
where $h(S)=f(S) \cap \mathfrak F$ for any simple group $S$. Thus,
if a formation
has an $\mathfrak L$-composition satellite, then it has an integrated
$\mathfrak L$-composition satellite.

\medskip
\textbf{Lemma 4.3.} \textsl{If  $\mathfrak F=CF_{\mathfrak L}(f)$, then
$\mathfrak F=CF_{\mathfrak L}(h)$  where $h$ is an integrated
 $\mathfrak L$-composition satellite such that 
$h(S)=\mathfrak F$ for every $S\in (\mathfrak L ^+)'$.  }

\medskip
\textbf{Proof.} We can  assume without loss of generality that $f$ is 
integrated. Let $\mathfrak H =CF_{\mathfrak L}(h)$  where $h(S)=f(S)$
if $S\in \mathfrak L ^+$, and $h(S)=\mathfrak F$ if $S\in (\mathfrak L ^+)'$.
Evidently, $\mathfrak F \subseteq \mathfrak H$. Assume that 
$\mathfrak H \not \subseteq \mathfrak F$, and choose a group $G$ of minimal 
order in $\mathfrak H \setminus \mathfrak F$. Then $L=G^{\mathfrak F}$ is a
unique minimal normal subgroup in $G$, and $L$ is not $f$-central.
Clearly, $c_G(L)=1$, and $C_G(L)=1$ if $L$ is non-abelian.  Let $A\in 
\mathrm{Com}(L)$. Applying Definition 4.4 and considering the cases 
$A\in \mathfrak L ^+$, $A\in \mathfrak L ^-$ and $A\in \mathfrak L '$,
we arrive at a contradiction. \qed

\medskip
\textbf{Theorem 4.4} (see \cite{ShSk2000}, Theorem 2). \textsl{  
Let $\mathfrak F$ be a non-empty formation, $\mathfrak L$ a class of simple
groups. The following statements are 
equivalent:}

\textsl{  
$(1)$ $\mathfrak F$ has an $\mathfrak L$-composition satellite; }

\textsl{  
$(2)$ $\mathfrak F$ has an ${\mathfrak L}^+$-composition satellite.}

\medskip

\textbf{Proof}.
 $ (1) \Rightarrow (2)$. Let $\mathfrak F=CF_{\mathfrak L}(f)$.
  Applying Lemma 4.3 we can suppose that $f$ is
integrated and $f(S)=\mathfrak F$ for every $S\in ({\mathfrak L}^+)'$. 
 Let  $\mathfrak H=CF_{{\mathfrak L}^+}(h)$ where $h$ is an 
${\mathfrak L}^+$-composition satellite such that $h(S)=f(S)$ if
$S\in {\mathfrak L}^+$, and $h(S)= \mathfrak F$ if
 $S\in {\mathfrak L}'\cup {\mathfrak L}^-=({{\mathfrak L}^+})'$. We will prove that
 $\mathfrak F=\mathfrak H$.

  If $G$ is a group of minimal order in $\mathfrak F\setminus \mathfrak H$, 
then $L=G^{\mathfrak H}$ is a unique minimal normal subgroup in $G$, and 
$L$ is not $h$-central. Clearly, $c_G(L)=1$, and $C_G(L)=1$ if $L$ is 
non-abelian. Applying Definition 4.4 we see that $L$ is $h$-central, a
contradiction. Thus  $\mathfrak F \subseteq \mathfrak H$.

Let $G$ be a group of minimal order in $\mathfrak H\setminus \mathfrak F$.
 Then $L=G^{\mathfrak F}$ is a unique minimal normal subgroup in $G$, and 
$L$ is not $f$-central. Clearly, $c_G(L)=1$, and $C_G(L)=1$ if $L$ is 
non-abelian. Applying again Definition 4.4 we see that $L$ is $f$-central, and
we arrive  at a
contradiction. Thus  $\mathfrak H \subseteq \mathfrak F$.

$ (2) \Rightarrow (1)$.  Let  $\mathfrak F=CF_{{\mathfrak L}^+}(f)$.
Applying Lemma 4.3 we can suppose that $f$ is
integrated and $f(({\mathfrak L}^+)')=\mathfrak F$. 
Let $h$ be an 
$\mathfrak L$-composition satellite such that $h(S)=f(S)$ if
$S\in {\mathfrak L}^+$, and
$h(S)=\mathfrak F$ if
$S\in ({\mathfrak L}^+)'$.  It is easy to
see that  $\mathfrak F=CF_{\mathfrak L}(h)$.    \qed

\medskip 
\textbf{Remark 4.4.}
 It follows from Theorem 4.4 that every non-empty formation 
$\mathfrak F$ with the property
 $\mathrm{Com}^+(\mathfrak F)\cap \mathfrak L =\varnothing$ has an 
$\mathfrak L$-composition satellite.

\medskip
\textbf{Remark 4.5.}
When $\mathfrak L = \mathfrak L ^+$ and $\omega = \pi (\mathfrak L)$, we usually
 use the term `$\omega$-composition satellite' and the notations 
$CF_{\omega}(f)$, $f(\omega ')$
in place of the term `$\mathfrak L$-composition satellite' and the 
notations $CF_{\mathfrak L}(f)$, $f(\mathfrak L ')$, respectively.

\medskip
\textbf{Theorem 4.5} {(see \cite{Shem97}, Theorems 3.1 and 3.2).\textsl{
Let $\mathfrak F$ be a non-empty formation, $\omega$ a set of primes. The following
statements  are pairwise equivalent: }

\textsl{$(1)$ $\mathfrak F$ is ${\mathfrak N}_{\omega}$-saturated; }

\textsl{$(2)$ $\mathfrak F$ is $\omega$-solubly saturated; }

\textsl{$(3)$  $\mathrm{cform}(\mathfrak F) \subseteq {\mathfrak N}_{{\omega}'} 
\mathfrak F$;}

\textsl{$(4)$ $\mathfrak F = CF_{\omega}(f)$ where  $f$ is a $\omega$-composition satellite
satisfying the following conditions:}

\textsl{$(i)$ $f(p)=\mathrm{Q}(G/C^p(G) \mid G\in \mathfrak F)$}
\textsl{if $p\in \omega$ and $C_p\in \mathrm{Com}(\mathfrak F)$;}

\textsl{$(ii)$ $f(p)=\varnothing$ if  $p\in \omega$ and
 $C_p\notin \mathrm{Com}(\mathfrak F)$; }

\textsl{$(iii)$ $f(S)=\mathfrak F$ if
 $S\in \mathfrak J \setminus \{C_p\mid p\in \omega \}$.}

\medskip
\textbf{Proof.} $(1)\Rightarrow (3)$. Set $\mathfrak H= \mathrm{cform}(\mathfrak F)$.
Fix $p\in \omega$. Since  $\mathfrak H \subseteq \mathfrak {NF}$, it is 
sufficient to show that $\mathfrak H \subseteq {\mathfrak N}_{p'} 
\mathfrak F$.
 Let $G$ be a group of minimal order in
$\mathfrak H \setminus {\mathfrak N}_{p'} \mathfrak F$. Clearly, $G$ is
monolithic and $L=\mathrm{Soc}(G)$ is the ${\mathfrak N}_{p'} \mathfrak F$-residual
of $G$. Since $\mathfrak F \subseteq {\mathfrak N}_{p'} \mathfrak F$, it follows that
$G^{\mathfrak F}\ge L$. Since $G\in \mathfrak H \subseteq \mathfrak {NF}$, we have
$G^{\mathfrak F} \in \mathfrak N$.  Since $G$ is monolithic and
 $G\notin {\mathfrak N}_{p'} \mathfrak F$, it follows that $G^{\mathfrak F}$ is
a $p$-group. From $G/L\in {\mathfrak N}_{p'} \mathfrak F$ it follows that 
$(G/L)^{\mathfrak F}= G^{\mathfrak F}/L$ is a $p'$-group. Therefore, 
$G^{\mathfrak F}=L=G^{{\mathfrak N}_{p'} \mathfrak F }$. By Lemma 4.2,
$\mathfrak H$ has a composition satellite $h$ such that 
$h(p)=\mathrm{Q}(A/C^p(A) \mid A\in \mathfrak F)$. Since $L$ is a 
$p$-group, we have $C_p\in \mathrm{Com}(G)$. Now from Lemma 4.1 it follows that 
$C_p\in \mathrm{Com}(\mathfrak F).$ Thus, applying Lemma 2.2, it follows that
${\mathfrak N}_p\subseteq \mathfrak F$. Since $h$ is a composition 
satellite of $\mathfrak H$, we have that $G/C_G(L)\in h(p)$. Therefore
$[L](G/C_G(L))\in \mathfrak H$, and $G/C_G(L)$ acts fixed-point-free on 
$L$. It follows that $G/C_G(L)\simeq T/N$, $T=A/C^p(A)$, $A\in\mathfrak F$.
If $C_p\notin \mathrm{Com}(A)$, then $A=C_p(A)$, $T=1$ and $G=L\in 
\mathfrak F$. Assume that $C_p\in \mathrm{Com}(A)$.  Since 
$G/C_G(L)\simeq T/N$, we can consider $L$ as an irreducible $\mathbb{F}_p(T/N)$-module.
Then $L$ becomes an irreducible $\mathbb{F}_pT$-module by inflation (see
\cite{DH}, p. 105). Since $T=A/C^p(A)$, we have by Lemma 2.5 that $[L]T\in \mathfrak F$.
By Lemma 2.6 it then follows that $[L](T/N)\in\mathfrak F$. From this and 
$T/N\simeq G/C_G(L)$ we deduce that $[L](G/C_G(L)\in \mathfrak F$. Hence,
by Lemma 2.4  it follows that $G\in \mathfrak F$. 

$(3)\Rightarrow (2)$. It is sufficient to consider only the case $\omega 
=\{p\}$. Let $H$ be a $p$-soluble normal subgroup in $G$,
   $L=O_p(H)\cap \Phi(H)$, and $G/L\in \mathfrak F$. 
We need to prove that $G\in \mathfrak F$. If $O_{p'}(H)\ne 1$, then
$LO_{p'}(H)/O_{p'}(H)\le \Phi(H/O_{p'}(H))$, and by induction we have
$G/O_{p'}(H)\in \mathfrak F$. From this and $G/L\in \mathfrak F$ it 
follows that $G\in \mathfrak F$. Assume that $O_{p'}(H)=1$.
By Lemma 4.2,
$\mathfrak H= \mathrm{cform}(\mathfrak F)$ has a composition satellite $h$ such that 
$h(p)=\mathrm{Q}(A/C^p(A) \mid A\in \mathfrak F)$.
 Let $t$ be a 
local satellite such that $t(p)=h(p)$ and $t(q)=\mathfrak E$ for every 
prime $q\ne p$.  Since $G/L\in \mathfrak F\subseteq \mathfrak H$ and $L=\Phi (H)$, $G$ acts 
$t$-hypercentrally on $H/\Phi (H)$. By Theorem 3.1, $G$ acts 
$t$-hypercentrally on $L=\Phi (H)$. But then $G$ acts 
$h$-hypercentrally on $L=\Phi (H)$, and 
we get $G\in \mathfrak H\subseteq {\mathfrak N}_{p'} \mathfrak F$. Thus,
$G^{\mathfrak F}\in {\mathfrak N}_{p'}\cap {\mathfrak N}_{p}=(1)$.   So,
$G\in \mathfrak F$, as required.

$(1)\Rightarrow (4)$.  Assume that $\mathfrak F$ is ${\mathfrak N}_{\omega}$-saturated.
Let $h$ be the minimal composition satellite of
 $\mathfrak H=\mathrm{cform}(\mathfrak F)$.  
Let $\mathfrak M =CF_{\omega}(f)$ where $f$ is an $\omega$-composition satellite
satisfying the following conditions:

 1) $f(p)=h(p)$ if $p\in \omega$;

2) $f(S)=\mathfrak F$ if $S\in \mathfrak J \setminus \{C_p\mid p\in \omega 
\}$.

 Inclusion $\mathfrak F \subseteq \mathfrak M$ 
is evident. Assume that the converse inclusion is false, and let $G$ be a 
group of minimal order in $\mathfrak M \setminus \mathfrak F$. Then 
$L=G^{\mathfrak F}$ is a unique minimal normal subgroup in $G$. If $L$ is
not an abelian $\omega$-group, it follows from $G\in \mathfrak M$ and 
$c_G(L)=1$ that
 $G\simeq G/c_G(L)\in \mathfrak F$.
 Therefore $L$ is
an $p$-group for some $p\in \omega$, and 
 we have $G/C^p(G)\in f(p)=h(p)$.
 Thus $G\in \mathfrak H$. Since $(1)\Rightarrow (3)$,  we get
 $G\in {\mathfrak N}_{p'}\mathfrak F$, and therefore 
$G^{\mathfrak F}\in {\mathfrak N}_{p'}\cap {\mathfrak N}_p = (1)$. So
$\mathfrak F = \mathfrak M$. We notice that by Lemma 4.2 we
have  $f(p)=h(p)=\varnothing$ if  $p\in \omega$ and
 $C_p\notin \mathrm{Com}(\mathfrak F)$.

$(4)\Rightarrow (1)$. Let $G/L\in \mathfrak F$ and $L= \Phi (O_p(G))$, $p\in \omega$.
 By Lemma 2.8, $C^p(G)/L=C^p(G/L)$.  Applying Proposition 4.1
to $G/L$, we have
$G/O_p(G)\simeq (G/L)/O_p((G/L))\in \mathfrak F$ and
$G/C^p(G)\simeq (G/L)/C^p(G/L)=(G/L)/C^p(G/L)\in f(p)$. But then by
Proposition 4.1 we get $G\in \mathfrak F$.  \qed
 
\medskip
\textbf{Corollary 4.5.1.} \textsl{If a non-empty formation $\mathfrak F$ is
$p$-solubly saturated and $C_p\in \mathrm{Com}(\mathfrak F)$,
 then $\mathfrak F$ has an $p$-composition satellite $f$ 
such that $f(p')=\mathfrak F$ and $f(p)=\mathrm{Q}(G/C^p(G)\mid G\in \mathfrak F)$.}

\medskip
\textbf{Corollary 4.5.2.} \textsl{If a non-empty formation $\mathfrak F$ is
solubly saturated, then $\mathfrak F= CF(f)$ where $f$ is a composition satellite 
satisfying the following conditions: }

\textsl{$(i)$  $f(p)=\mathrm{Q}(G/C^p(G)\mid G\in \mathfrak F)$ if $p\in\omega$ 
and $C_p\in \mathrm{Com}(\mathfrak F)$; }

\textsl{$(ii)$ $f(S)=\mathfrak F$ for every $S\in \mathrm{Com^-}(\mathfrak 
F)$;}

\textsl{$(iii)$ $f(S)=\varnothing$ for every $S\in \mathfrak J\setminus \mathrm{Com}(\mathfrak F).$ }
 
\medskip
\textbf{Theorem 4.6} (see \cite{Shem97}, Theorem 3.1(b)).
\textsl{Let $\mathfrak F$ be a non-empty $\omega$-saturated formation,
and $h$ be the minimal composition satellite of $\mathrm{cform}(\mathfrak F)$.
Then $\mathfrak F=LF_{\omega}(f)$ where $f$ is an $\omega$-local
satellite such that $f(p)=h(p)$ for every $p\in \omega$.      }  

\medskip
\textbf{Proof.} We may suppose without loss of generality that 
$\omega \subseteq \pi (\mathfrak F)$. By Lemma 4.2,
 $h(S) = \mathrm{Q}(H/C^S(H) \mid H\in \mathfrak F)$ if 
$S\in \mathrm{Com}(\mathfrak F)$, and $h(S)=\varnothing$
 if $S\in \mathfrak J \setminus \mathrm{Com}(\mathfrak F)$. 

Let $p$ be a prime in $\omega$, and $S$ be a non-abelian $pd$-group
 in $\mathrm{Com}(\mathfrak F)$.  We will now prove that
  $h(S)\subseteq h(p)$. Consider $R=H/C^S(H)$, $H\in \mathfrak F$. 
By Lemma 2.1, $C^S(H)$ is the largest normal subgroup not having 
composition factors isomorphic to $S$. Clearly, $O_{p',p}(R)=1$. 
Let $A_p(R)$ be the $p$-Frattini module, i.~e., the kernel of the 
universal Frattini, $p$-elementary $R$-extension:
\[
 1\rightarrow A_p(R) \xrightarrow {\mu} E \xrightarrow {\varepsilon}
 R\rightarrow 1.
\]
Here $E/A_p(R)\simeq R$, and $A_p(R)$ is an elementary abelian $p$-group
 contained in $\Phi (E)$. Let $N_1,\dots, N_t$ be all minimal normal
subgroups in $E$ contained in   $A_p(R)$. Since $\mathfrak F$
 is $p$-saturated, we have
$E\in \mathfrak F \subseteq \mathrm{cform}(\mathfrak F)$, and therefore
$E/\cap _i C_E(N_i)\in h(p)$. Since $N_1,\dots, N_t$ are simple submodules
of the $\mathbb{F}_pR$-module $A_p(R)$, it follows that
$R/\mathrm{Ker}(R \;on\;(N_1\dots N_t))\in h(p)$. By theorem 
of Griess and P.~Schmid, 
$\mathrm{Ker}(R \; on\; (N_1\dots N_t))=O_{p',p}(R)$
 (see \cite{GS} or \cite{DH}, p. 833). Since $O_{p',p}(R)=1$, it follows that
$R\in h(p)$. Thus, $h(S) = \mathrm{Q}(H/C^S(H) \mid H\in 
\mathfrak F)\subseteq h(p)$ if
 $S \in \mathrm{Com}(\mathfrak F)$ and $p\in \omega \cap \pi (S)$. 

  Let $f$ be an $\omega$-local
satellite such that $f(p)=h(p)$ if $p\in \omega$, and $f(\omega ')=\mathfrak F$ 
if ${\omega}'\ne \varnothing$. We will prove now 
that $\mathfrak F = LF_{\omega}(f)$. 

Let $G$ be a group of minimal order in $\mathfrak F \setminus LF_{\omega}(f)$.
Then $L=G^{LF_{\omega}(f)}$ is a unique minimal normal subgroup in $G$, and $L$
 is not $f$-central in $G$. If $L$ is an $\omega '$-group, then $c_G(L)=1$ 
and $G\simeq G/c_G(L)\in f(\omega ')=\mathfrak F$. If $L$ is a non-abelian
$pd$-group for some $p\in \omega$ and $S\in \mathrm{Com}(L)$, then $C_G(L)=1$
 and we have
$G\simeq G/C_G(L)\in h(S)\subseteq h(p) \subseteq \mathfrak F$.
 Assume that $L$ is a $p$-group, $p\in \omega$.
Since $L$ is not $f$-central, $L\not \subseteq Z(G)$. By Lemma 2.1 we have
$C^p(G)=1$. So $G\in h(p) = \mathrm{Q}(H/C^p(H) \mid H\in \mathfrak F)$, i.~e.,
$L$ is $f$-central, a contradiction. Thus, $\mathfrak F \subseteq LF_{\omega}(f)$.                        
     
Let $G$ be a group of minimal order in $LF_{\omega}(f)\setminus \mathfrak F$. Then 
 $L=G^{\mathfrak F}$ is a unique minimal normal subgroup in $G$.
Clearly, $c_G(L)=1$, and $C_G(L)=1$ if  $L$ is non-abelian. Hence,
it follows from  $G\in LF_{\omega}(f)$ that if $L$
is an $\omega '$-group, then $G\in f(\omega ')= \mathfrak F$, and if
$L$ is a non-abelian $pd$-group for some $p\in \omega$, then  $G\in f(p)=h(p)\subseteq 
\mathfrak F$, and we get a contradiction. 
 Assume that 
$L$ is a $p$-group, $p\in \omega$. Evidently, $L$ is not contained in 
$\Phi (G)$ (recall that $\mathfrak F$ is $p$-saturated). By Lemma 2.7, 
$L=C_G(L)$.  Since $L$ is $f$-central, we obtain that $G=[L]T$ where $T\in 
f(p)$.  Therefore, $T\simeq R/K$ where $R=H/C^p(H)$ for some $H\in 
\mathfrak F$. Now we can consider $L$ as an irreducible ${\mathbb 
F}_pR$-module by inflation (see \cite{DH}, p. 105). By Lemma 2.5 we have
$[L]R\in \mathfrak F$. Since $K$ acts identically on $L$, it follows
from  Lemma 2.6 that  $[L](R/K)\simeq LT=G\in 
\mathfrak F$, and we again arrive at a contradiction.
So $LF_{\omega}(f)= \mathfrak F$. \qed

\medskip
\textbf{Corollary 4.6.1.} \textsl{If a non-empty formation $\mathfrak F$ is
$\omega$-saturated, then $\mathfrak F$ has an $\omega$-local satellite $f$ 
such that  $f(p)=\mathrm{Q}(G/C^p(G)\mid G\in \mathfrak F)$ if 
$p\in \omega \cap \pi (\mathfrak F)$, $f(p)=\varnothing$ if
 $p\in \omega \setminus \pi (\mathfrak F)$, and $f({\omega}')=\mathfrak F$ if 
${\omega}'\ne \varnothing$.}

\medskip
\textbf{Corollary 4.6.2.} \textsl{If a non-empty formation $\mathfrak F$ is
saturated, then $\mathfrak F= LF(f)$ where $f$ is a local satellite  
such that  $f(p)=\mathrm{Q}(G/C^p(G)\mid G\in \mathfrak F)$ for every 
$p\in \pi (\mathfrak F)$, and $f(p)=\varnothing$ for every prime
 $p\notin \pi (\mathfrak F)$}. 

\bigskip
\textbf{4. $\mathfrak X$-local formations }

\bigskip
In 1985 F\"{o}rster \cite{For85} introduced the concept `$\mathfrak X$-local 
formation' in order to obtain a common extension of Theorem 3.2 and 4.1.

\medskip
\textbf{Definition 5.1. }
Let $\mathfrak X$ be a class of simple groups such that 
  ${\mathrm {Char}}(\mathfrak X) = \pi (\mathfrak X)$. Consider a map
\[
   f: \pi ( \mathfrak X)\cup \mathfrak X' \longrightarrow \{\text{formations}\}
\]
which does not distinguish between any two non-identity isomorphic groups. Denote 
through $LF_{\mathfrak X}(f)$ the class of all groups  $G$ satisfying
the following conditions:

(i) if $H/K$ is a chief ${\mathrm E}\mathfrak X$-factor of a group $G$, then
 $G/C_G (H/K)$ belongs to $f(p)$ for any $p \in \pi (H/K)$;

(ii) if $G/L$ is a  monolithic quotient of $G$
 and $\mathrm {Soc}(G/L) \in {\mathrm E}(\mathfrak X ')$, then
$G/L \in f(S)$ where $S \in \mathrm {Com}(\mathrm {Soc}(G/L))$.

The class $LF_{\mathfrak X}(f)$ is a formation; it is called an
$\mathfrak X${\it -local formation}. 

\medskip
$\mathfrak X$-local formations were investigated in \cite{BCE05,
BB91,BBC,BShem,BE}. In \cite{Shem10} it was 
proved with help of some lemmas in \cite{BShem} that every
 $\mathfrak X$-local 
formation has a ${\mathfrak X}^+$-composition satellite. Now we give a direct
proof of that fact.

\medskip
\textbf{Theorem 5.1.} \textsl{Let $\mathfrak F$ be a non-empty formation,
 $\mathfrak X$  a class of simple groups such that 
  ${\mathrm {Char}}(\mathfrak X) = \pi (\mathfrak X)$. Let $\mathfrak L$ be a
class of simple groups such that $\mathfrak L ^+=\mathfrak X ^+$. }

\textsl{$(1)$ If  $\mathfrak F$ is an $\mathfrak X$-local formation, then
$\mathfrak F$ has an $\mathfrak L$-composition satellite.}

\textsl{$(2)$ If  $\mathfrak F$ has an $\mathfrak L$-composition 
satellite, then $\mathfrak F$ is an $\mathfrak X ^+$-local formation.}
 
\medskip
\textbf{Proof.}
Set  $\omega ={\mathrm {Char}}(\mathfrak X)$. Evidently,
$\mathfrak L ^- \cup \mathfrak L ' = \mathfrak X ^- \cup \mathfrak X '=
 (\mathfrak L ^+)'=(\mathfrak X ^+)'$. 

 (1) Let $\mathfrak F$ be a $\mathfrak 
X$-local formation, $\mathfrak F = LF_{\mathfrak X}(f)$. Consider an
$\mathfrak L$-composition satellite $h$ such that $h(p)=f(p)\cap \mathfrak F$ if
$p\in \omega$, and $h(S)=\mathfrak F$ if $S\in \mathfrak L ^- \cup 
\mathfrak L '$. We will prove that $\mathfrak F = CF_{\mathfrak L}(h)$.

Suppose that $\mathfrak F \not \subseteq CF_{\mathfrak L}(h)$. Let $G$ be 
a group of minimal order in $\mathfrak F \setminus CF_{\mathfrak L}(h)$. 
Then $G$ is monolithic, and $G/M \in CF_{\mathfrak L}(h)$ where $M$ is the 
socle of $G$. Clearly, $M$ is the $CF_{\mathfrak L}(h)$-residual of $G$, 
and every chief factor between $G$ and $L$ is $h$-central. Assume that $M$ is an
$\mathrm{E}(\mathfrak L ^- \cup \mathfrak L ')$-group. Since $G\in \mathfrak F$,
we have that $G\in  h(S)$ where $S\in \mathrm{Com}(M)$. Since 
$c_G(L)=1$,  we have that $M$ is $h$-central in $G$, and so
 $G\in CF_{\mathfrak L}(h)$. Assume now that  $M$ is a $p$-group, $p\in \omega$.
Since $G\in \mathfrak F$, we have $G/C_G(M)\in f(p)\cap 
\mathfrak F = h(p)$, i.~e.,  $M$ is $h$-central. We 
see that $G\in CF_{\mathfrak L}(h)$, a contradiction. Thus, 
$\mathfrak F  \subseteq CF_{\mathfrak L}(h)$.

 Suppose now that $CF_{\mathfrak L}(h)\not \subseteq \mathfrak F$. Choose a 
group $G$ of minimal order in $CF_{\mathfrak L}(h)\setminus \mathfrak 
F$. Then $G$ is monolithic, and $G/M\in \mathfrak F$ where $M=G^{\mathfrak 
F}$ is the socle of $G$. Assume that $M$ is an  $\mathrm{E}(\mathfrak L ^- \cup 
\mathfrak L ')$-group. Then from $c_G(L)=1$ and $h$-centrality of $L$ it 
follows that $G/c_G(M)\simeq G\in \mathfrak F$. Assume that $M$ is a $p$-group,
 $p\in \omega$. Then
 \[
 G/C_G(M)\in h(p)= \mathfrak F \cap f(p)\subseteq f(p).
 \]
 We see  that all the chief factors and all the quotients of $G$
  satisfies conditions (1) and (2) of Definition 5.1. So, $G\in \mathfrak 
F$, a contradiction. Thus, $\mathfrak F = CF_{\mathfrak L}(h)$.

(2) Let $\mathfrak F$ be a formation having an $\mathfrak 
L$-composition satellite.  By Lemma 4.3,  $\mathfrak F = CF_{\mathfrak 
L}(f)$ where $f$ is an $\mathfrak L$-composition satellite such that
$f(S)=\mathfrak F$ for every $S\in \mathfrak L ^- \cup \mathfrak L '$.
Consider an $\mathfrak X ^+$-local formation  $\mathfrak H = LF_{\mathfrak 
X ^+}(h)$ where $h(p)=f(p)$ for any $p\in \omega$, and $h(S)=\mathfrak F$ for
every $S\in (\mathfrak X ^+)'$.  It easy to check that 
$\mathfrak F = \mathfrak H$. \qed

\end{document}